\documentclass[11pt]{article}
\usepackage{amsmath}
\usepackage{amsmath,amssymb,amsfonts,amsthm,fancyhdr}
\usepackage{epsfig,graphicx,picins,picinpar,subfigure}
\usepackage{pstricks}
\usepackage{fancyvrb}
\usepackage[numbers,sort&compress]{natbib}
\begin{document}


\title{\textbf{Notes on Dickson's Conjecture}}
\author{\emph{Shaohua Zhang}}
\date{{\small  School of Mathematics, Shandong University,
Jinan,  Shandong, 250100, China \\E-mail:
shaohuazhang@mail.sdu.edu.cn}}
\maketitle

\begin{abstract}
In 1904, Dickson [5] stated a very important conjecture. Now people
call it Dickson's conjecture. In 1958, Schinzel and Sierpinski [14]
generalized Dickson's conjecture to the higher order integral
polynomial case. However, they did not generalize Dickson's
conjecture to the multivariable case. In 2006, Green and Tao [13]
considered Dickson's conjecture in the multivariable case and gave
directly a generalized Hardy-Littlewood estimation. But, the precise
Dickson's conjecture in the multivariable case does not seem to have
been formulated. In this paper, based on the idea in [15], we will
try to complement this and give an equivalent form of Dickson's
Conjecture, furthermore, generalize it to the multivariable case or
a system of affine-linear forms on $N^k$ . We also give some remarks
and evidences on conjectures in [15]. Finally, in Appendix, we
briefly introduce the basic theory that several multivariable
integral polynomials represent simultaneously prime numbers for
infinitely many integral points.

\vspace{3mm}\textbf{Keywords:} Chinese Remainder Theorem,
Dirichlet's theorem, Dickson's conjecture, Green-Tao theorem,
affine-linear form

\vspace{3mm}\textbf{2000 MR  Subject Classification:}\quad 11A41,
11A99
\end{abstract}


\section{Introduction}
\setcounter{section}{1}\setcounter{equation}{0}
The question of existence of infinitely many prime values of
polynomials $f(x)$ with integral coefficients has been one of the
most important topics in Number Theory.  Euclid [1] proved firstly
that $f(x)=x$ represents infinitely many primes. In 1837, Dirichlet
[2] showed that $f(x)=a+bx$ takes infinitely many primes, where $a$
and $b$ are integers satisfying $(a,b)=1$, and either $a\neq 0,
b>0$, or $a=0, b=1$. In 1857, Bouniakowsky [3] considered the case
of nonlinear polynomials and conjectured that if $f(x)$ is an
irreducible polynomial with integral coefficients, positive leading
coefficient and degree at least 2, and there does not exist any
integer $n>1$ dividing all the values $f(k)$ for every integer $k$,
then $f(x)$ is prime for an infinite number of integers $x$.
Unfortunately, as far, his conjecture even the simplest case
$f(x)=x^2+1$ [4] is still open. In a somewhat different direction,
by generalizing Dirichlet's theorem and concerning the simultaneous
values of several linear polynomials, Dickson [5] stated the
following conjecture in 1904:

\vspace{3mm}\noindent {\bf Dickson's conjecture:~~}%
Let $s\geq1$, $f_i(x)=a_i+b_ix$ with $a_i$ and $b_i$ integers,
$b_i\geq1$ (for $i=1,...,s$ ). If there does not exist any integer
$n>1$ dividing all the products $\prod_{i=1}^{i=s}f_i(k)$, for every
integer $k$, then there exist infinitely many natural numbers $m$
such that all numbers $f_1(m),...,f_s(m)$ are primes.

\vspace{3mm} Dickson's conjecture implies many important results [6]
such as:

1, there exist infinitely many composite Mersenne numbers.

2, there exist infinitely many pairs of twin primes.

3, there exist infinitely many Carmichael numbers.

4, Artin's conjecture is true.

5, a conjecture of Hardy and Littlewood is true.

6, there exist infinitely many Sophie Germain primes or safe primes.

7, van der Corput's theorem [7] which states that there are
infinitely many triples of primes in arithmetic progression.

8, Balog's theorem [10] which states that for any $m>1$, there are
$m$ distinct primes $p_1,...,p_m$  such that all of the averages
$\frac{p_i+p_j}{2}$ are primes.

9, Green-Tao theorem [8] which states that the sequence of prime
numbers contains arbitrarily long arithmetic progressions.

10, and so on.

\vspace{3mm}Historically, many mathematicians were interesting in
Dickson's conjecture and obtained great results. For example, Hardy
and Littlewood [9], Balog [10, 11], Heath-Brown [12], Green and Tao
[13], and et al. gave important consideration and profound analysis
on Dickson's conjecture and its special cases. In [13], Green and
Tao further generalized the conjecture of Hardy and Littlewood in
[9].

\vspace{3mm}In 1958, by studying the consequences of Bouniakowsky's
conjecture and Dickson's conjecture, A. Schinzel and W. Sierpinski
[14] got the following:

\vspace{3mm}\noindent {\bf Schinzel-Sierpinski conjecture (H hypothesis):~~}%
Let  $s\geq1$, and let $f_1(x),...,f_s(x)$ be irreducible
polynomials with integral coefficients and positive leading
coefficient. If there does not exist any integer $n>1$ dividing all
the products $\prod_{i=1}^{i=s}f_i(k)$, for every integer $k$, then
there exist infinitely many natural numbers $m$ such that all
numbers $f_1(m),...,f_s(m)$ are primes.

\vspace{3mm} However, Schinzel and Sierpinski did not generalize
Dickson's conjecture to the multivariable case. In [13], Green and
Tao considered Dickson's conjecture in the multivariable case. But,
the precise conjecture does not seem to have been formulated. In
this paper, based on the idea in [15], we will try to complement
this and give an equivalent form of Dickson's conjecture, and,
generalize it to the multivariable case or a system of affine-linear
forms on $N^k$. In this paper, we always restrict that a
$k$-variables integral polynomial is a map from $N^k$ to $Z$, where
$k\geq 1$. We will also give some remarks and evidences on
conjectures in [15]. We obtained the following:

\vspace{3mm}\noindent {\bf The equivalent form of Dickson's conjecture:~~}%
Let $f_1(x),...,f_s(x)$ be $s$ linear polynomials with integral
coefficients, if there is a positive integer $c$ such that for every
positive integer $m\geq c$, there exists a positive integer $y$ such
that $f_1(y)>1,...,f_s(y)>1$ are all in $Z_m^* =\{x|1\leq x <m,
(x,m)=1\}$, then $f_1(x),...,f_s(x)$ represent simultaneously primes
for infinitely many positive integers $x$.

\vspace{3mm}Namely, the sufficient and necessary condition that $s$
linear polynomials  $f_1(x),...,f_s(x)$ with integral coefficients
 represent infinitely many primes is of that there is a positive integer $c$ such that for every
integer $m\geq c$, there exists a positive integer $y$ such that
$f_1(y)>1,...,f_s(y)>1$ are all in $Z_m^*$. This  sufficient and
necessary condition implies that there does not exist any integer
$n>1$ dividing all the products $\prod_{i=1}^{i=s}f_i(k)$, for every
integer $k$, and, for every $1\leq i\leq s$, the leading coefficient
of $f_i(x)$ is positive. Of course, it also implies the non-trivial
case that all the polynomials $f_i(x)$ are non-constant, and, no two
polynomials are rational multiples of each other.

\vspace{3mm}\noindent {\bf The generalization of Dickson's conjecture:~~}%
Let $s, k \in N$ and let $f_1(x_1,...,x_k),...,f_s(x_1,...,x_k)$  be
multivariable polynomials of degree 1 with  integral coefficients,
if there is a positive integer $c$ such that for every positive
integer $m\geq c$, there exists an integral point $(y_1,...,y_k)$
such that $f_1(y_1,...,y_k)>1,...,f_s(y_1,...,y_k)>1$ are all in
$Z_m^*$, then $f_1, ... , f_s$ represent simultaneously primes for
infinitely many integral points $(x_1,...,x_k)$.

\section{The proof of equivalent form}
\vspace{3mm}\noindent{\bf  The proof of case $s=1$:~~}%
We firstly prove that the equivalent form of Dickson's conjecture
holds when $s=1$. Namely, the sufficient and necessary condition
that the linear polynomial  $f(x)$ with integral coefficients
 represents infinitely many primes is of that there is a positive integer $c$ such that for every
integer $m\geq c$, there exists a positive integer $y$ such that
$f(y)>1$ is in $Z_m^*$. Let $f(x)=a+bx$ be a linear polynomial with
integral coefficients. Obviously, if there is a positive integer $c$
such that for every integer $m\geq c$, there exists a positive
integer $y$ such that $f(y)>1$ is in $Z_m^*$, then either $a=0, b=1$
or $a\neq 0, b>0$ with $(a,b)=1$. By Euclid's second theorem and
Dirichlet's theorem, $f(x)=a+bx$ represents infinitely many primes.

\vspace{3mm} On the other hand, if $f(x)=a+bx$ represents infinitely
many primes, then either $a=0, b=1$, or $a\neq 0, b>0$ with
$(a,b)=1$. When $a=0, b=1$, let $c=3$, then for every integer $m\geq
c$, there exists a positive integer $y$ e.g. $y=m-1$ such that
$f(y)=y=m-1>1$ is in $Z_m^*$. When $a\neq 0, b>0$ with $(a,b)=1$, in
[15], by using the result of de la Vall\'{e}e Poussin, which states
that the number of primes of the form $a+bn$ not exceeding a large
number $x$ is asymptotic to $x/{\varphi(b)\log x}$ as
$x\rightarrow\infty$, where $\varphi(.)$ is Euler's function, we
proved that there is a positive integer $c$ such that for every
integer $m\geq c$, there exists a positive integer $y$ such that
$f(y)>1$ is in $Z_m^*$. Next, we will give another proof which is
very simple and elementary.

\vspace{3mm} Note that $f(x)=a+bx$ represents infinitely many primes
when $a\neq 0, b>0$ with $(a,b)=1$. Therefore, there are positive
integers $y,z$ such that $p=a+by>q=a+bz>b(|a|+k)$, where $p,q$ are
primes and $k$ is the least positive integer such that $k>1$ and
$(a,k)=1$. Surely, $k$ is prime, too. Now, we prove that for every
integer $m\geq c=pq$, there exists a positive integer $x$ such that
$f(x)>1$ is in $Z_m^*$. Assume that $m\geq c=pq$. If $(m,p)=1$, then
we can choose $x=y$ such that $f(x)>1$ is in $Z_m^*$. Similarly, if
$(m,q)=1$, we can choose $x=z$. Hence, we only consider the case
$m=pqt$, where $t\in N$. Since $t\in N$, $a\in Z$ and $a\neq 0$.
Hence, $t$ has a positive divisor which is co-prime to $a$. For
example, $t$ has a positive divisor 1 and $(a,1)=1$. Let $d$ be the
greatest positive divisor of $t$ such that $(a,d)=1$. Then, either
$(bpd+a,pqt)=1$ or $(b(|a|+k)pd+a, pqt)=1$. We will prove this key
fact.

\vspace{3mm} We write $t=dr$ with $r\in N$. If $(a,t)=1$, then
$d=t,r=1$ because $d$ is the greatest positive divisor of $t$ such
that $(a,d)=1$. Thus, $(bpd+a,pt)=(bpt+a,pt)=1$ and $(b(|a|+k)pd+a,
pt)=1$ since $p>|a|$ and $(p,a)=1$. If $(a,t)\neq1$, then $r>1$ and
any prime divisor of $r$ divides $a$. But,
$(bpd,a)=(b(|a|+k)pd,a)=1$. Therefore, $(bpd+a,pt)=(bpd+a,pdr)=1$
and $(b(|a|+k)pd+a, pt)=1$. Note that $q=a+bz>b(|a|+k)$. So, $q$ can
not divide simultaneously $bpd+a$ and $b(|a|+k)pd+a$. So, either
$(bpd+a,pqt)=1$ or $(b(|a|+k)pd+a, pqt)=1$. On the other hand,
clearly, $1<bpd+a<pqt$ and $1<b(|a|+k)pd+a<pqt$. Thus, we can choose
a number $x$ ($x=pd$ or $x=(|a|+k)pd$) such that $f(x)>1$ is in
$Z_m^*$.

\vspace{3mm} In the proof above, we used Dirichlet's theorem.
However, one can give the third proof without Dirichlet's theorem.
Namely, without Dirichlet's theorem, one might prove directly that
if $f(x)=a+bx$ with $a\neq 0, b>0$ and $(a,b)=1$, then there is a
positive integer $c$ such that for every integer $m\geq c$, there
exists a positive integer $y$ such that $f(y)>1$ is in $Z_m^*$.

\vspace{3mm} Since $a\neq 0, b>0$, hence there is a positive integer
$u$ such that $f(u)=a+bu>1$, $(u,a)=1$, and any prime divisor of
$a+bu$ is greater than $b(|a|+k)$, where $k$ is the least positive
integer such that $k>1$ and $(a,k)=1$. For example, one can choose
$u=\prod_{p\leq 2b(|a|+k), (a,p)=1} p$, where $p$ is prime. Let $P$
be the set of all prime divisors of $a+bu$. On the other hand,
clearly, there is a positive integer $v$ such that $f(v)=a+bv>1$,
$(a+bu,a+bv)=1$ and any prime divisor of $a+bv$ is greater than
$b(|a|+k)$. For instance, one might choose $v=\prod_{p\in P} p
\prod_{p\leq 2b(|a|+k), (a,p)=1} p$, where $p$ is prime. Let $Q$ be
the set of all prime divisors of $a+bv$. We claim that for every
integer $m\geq c=(a+bu)\times (a+bv)$, there exists a positive
integer $y$ such that $f(y)>1$ is in $Z_m^*$. By the aforementioned
idea, it is enough to consider the case $m=pqt$, where $p\in P$,
$q\in Q$ and $t\in N$. The details of proof are left as an exercise.

\vspace{3mm}\noindent{\bf Remark 1:~~}%
It is worthwhile pointing out that the problem of determining the
low bound of constant $c$ is interesting.  For any integers $a,b$
satisfying $a\neq 0, b>0$ and $(a,b)=1$, let $c$ be the least
positive integer such that for every integer $m\geq c$, there exists
a positive integer $x$ such that $f(x)=a+bx>1$ is in $Z_m^*$. We
write $(r-1)!<c\leq r!$, where $r$ is a positive integer. By the
idea in [15], we conjecture that if $f(x)=a+bx>1$ in $Z_{r!}^*$ is
the least positive integer of the form $a+bx$, then it always is
prime. And if this holds, then it leads to a new proof on
Dirichlet's theorem. Moreover, as a special case of a conjecture in
[15],  the following conjecture is consistent with Dirichlet's
theorem.

\vspace{3mm}\noindent{\bf  Conjecture 1:~~}%
With the notation above, if $n\geq r$, then there always is a prime
of the form $a+bx$ in $Z_{n!}^*$.

\vspace{3mm}\noindent{\bf  The proof of case $s>1$:~~}%
Clearly, the condition that there is a positive integer $c$ such
that for every integer $m\geq c$, there exists a positive integer
$y$ such that $f_1(y)>1,...,f_s(y)>1$ are all in $Z_m^*$ implies
that there does not exist any integer $n>1$ dividing all the
products $\prod_{i=1}^{i=s}f_i(k)$, for every integer $k$, and, for
every $1\leq i\leq s$, the leading coefficient of $f_i(x)$ is
positive. Therefore, it is enough to prove that the latter implies
the former. For this goal, we write $f_i(x)=a_i+b_ix$ with $a_i$ and
$b_i$ integers, $b_i\geq1$ (for $i=1,...,s$ ). Since there does not
exist any integer $n>1$ dividing all the products
$\prod_{i=1}^{i=s}f_i(k)$, for every integer $k$, hence, $(a_i,
b_i)=1$ for every $1\leq i\leq s$. Without loss of generality, we
only need consider the two cases as follows:

\vspace{3mm}\noindent{\bf  Case 1:~~}%
Let $s>1$, $f_i(x)=a_i+b_ix$ with $a_i\neq 0$ and $b_i\geq 1$
integers for $i=1,...,s$ . If there does not exist any integer $n>1$
dividing all the products $\prod_{i=1}^{i=s}f_i(k)$, for every
integer $k$, then, there is a positive integer $c$ such that for
every positive integer $m\geq c$, there exists a positive integer
$y$ such that $f_1(y)>1,...,f_s(y)>1$ are all in $Z_m^*$.

\vspace{3mm}\noindent{\bf  Case 2:~~}%
Let $s>1$, $f_1(x)=x$, $f_i(x)=a_i+b_ix$ with $a_i\neq 0$ and
$b_i\geq 1$ integers for $i=2,...,s$. If there does not exist any
integer $n>1$ dividing all the products $\prod_{i=1}^{i=s}f_i(k)$,
for every integer $k$, then, there is a positive integer $c$ such
that for every positive integer $m\geq c$, there exists a positive
integer $y$ such that $f_1(y)>1,...,f_s(y)>1$ are all in $Z_m^*$.

\vspace{3mm}\noindent{\bf  Lemma 1:~~}%
Let $s\geq 1$, $f_i(x)=a_i+b_ix$ with $a_i\neq 0$, $b_i\geq 1$ and
$(a_i,b_i)=1$ for $i=1,...,s$ . If there does not exist any integer
$n>1$ dividing all the products $\prod_{i=1}^{i=s}f_i(k)$, for every
integer $k$, then there exists a positive integer $x$ such that the
least prime divisor of $\prod_{i=1}^{i=s}f_i(x)$ is greater than any
given positive integer $C$.

\vspace{3mm}\noindent{\bf  The proof of Lemma 1:~~}%
Easy.  When $C<2$, it is clear. Using the method in [15], we write
$C!=\prod_{i=1}^{i=r}{p_i}^{e_i}$ when  $C\geq 2$. Noticed that
there does not exist any integer $n>1$ dividing all the products
$\prod_{i=1}^{i=s}f_i(k)$, for every integer $k$. So, there exists a
positive integer
 $a_j$ such that $\gcd (\prod_{i=1}^{i=s}f_i(a_j),{p_j}^{e_j})=1 $  for $1\leq j\leq r$.
 By Chinese Remainder Theorem, there exists a positive integer $x$ such that $x\equiv a_j
(\mod {p_j}^{e_j})$. Note that $f_i(x)$ is a  polynomial with
integral coefficients. Hence, $\prod_{i=1}^{i=s}f_i(x)\equiv
\prod_{i=1}^{i=s}f_i(a_j) (\mod {p_j}^{e_j})$ and $\gcd
(\prod_{i=1}^{i=s}f_i(x),C!)=1$. It shows immediately that Lemma 1
holds.

\vspace{3mm}\noindent{\bf  Lemma 2:~~}%
Let $s\geq 1$,  $f_i(x)=a_i+b_ix$ with $a_i\neq 0$, $b_i\geq 1$ and
$(a_i,b_i)=1$ for $i=1,...,s$, and let $r$ be a positive integer
satisfying $(r,\prod_{i=1}^{i=s}a_i)=1$. If there does not exist any
integer $n>1$ dividing all the products $\prod_{i=1}^{i=s}f_i(k)$,
for every integer $k$, then for any positive integer $m$, there
exists a positive integer $x$ such that
$(\prod_{i=1}^{i=s}(a_i+r\times b_ix),m)=1$.

\vspace{3mm}\noindent{\bf  The proof of Lemma 2:~~}%
Clearly, by Chinese Remainder Theorem, it is enough to prove that
for any prime $p$, there exists a positive integer $x$ such that
$(\prod_{i=1}^{i=s}(a_i+r\times b_ix),p)=1$. If $p|r$, due to
$(r,\prod_{i=1}^{i=s}a_i)=1$, we can choose $x=1$ such that
$(\prod_{i=1}^{i=s}(a_i+r\times b_i),p)=1$. Now, we consider the
case that $p$ does not divide $r$. Using the known condition that
there does not exist any integer $n>1$ dividing all the products
$\prod_{i=1}^{i=s}f_i(k)$, for every integer $k$, we deduce easily
that there exists a positive integer $y$ such that
$(\prod_{i=1}^{i=s}(a_i+b_iy),p)=1$. Note that there are positive
integer $t,x$ such that $y+pt=rx$. So,
$(\prod_{i=1}^{i=s}(a_i+r\times b_ix),p)=1$. And Lemma 2 holds.

\vspace{3mm}\noindent{\bf  Lemma 3:~~}%
Let $s\geq 1$,  $f_i(x)=a_i+b_ix$ with $a_i\neq 0$, $b_i\geq 1$ and
$(a_i,b_i)=1$ for $i=1,...,s$, and let $r>1$ be a positive integer
satisfying $(r,\prod_{i=1}^{i=s}a_i)=1$.  If there does not exist
any integer $n>1$ dividing all the products
$\prod_{i=1}^{i=s}f_i(k)$, for every integer $k$, then for any
positive integer $l$, there is a positive integer $e>1$ depending on
$l,s,a_1,...,a_s$, such that for any positive integer $m$, there
exists a positive integer $x$ satisfying
$(\prod_{j=0}^{j=l}\prod_{i=1}^{i=s}(a_i+r^{e^j}\times b_ix),m)=1$.

\vspace{3mm}\noindent{\bf  The proof of Lemma 3:~~}%
For any positive integer $l$, we choose $f=\prod_{p\leq
s^{l+1}\prod_{i=1}^{i=s}|a_i|+1,(p,\prod_{i=1}^{i=s}|a_i|)=1}p$,
where $p$ represents a prime. Let $e=\varphi (f)+1$. Then for any
prime $q$ satisfying $q\leq s^{l+1}$, $r^{e^j}\equiv r (\mod q)$
holds for $j=0,...,l$. By the known condition and Lemma 2, we know
that there exists a positive integer $x$ such that
$(\prod_{i=1}^{i=s}(a_i+r\times b_ix),q)=1$. So,
$(\prod_{j=0}^{j=l}\prod_{i=1}^{i=s}(a_i+r^{e^j}\times b_ix),q)=1$.
If $q>  s^{l+1}$, clearly, there exists a positive integer $x$
satisfying $(\prod_{j=0}^{j=l}\prod_{i=1}^{i=s}(a_i+r^{e^j}\times
b_ix),q)=1$ since the degree of polynomial
$(\prod_{j=0}^{j=l}\prod_{i=1}^{i=s}(a_i+r^{e^j}\times b_ix)(\mod
q)$ is at most $ s^{l+1}$. By Chinese Remainder Theorem, one can
further prove the Lemma 3 holds for $s>1$. While $s=1$, Lemma 3 is
trivial by picking $e=\varphi (\prod_{p\leq
(l+1)|a_1|+1,(p,|a_i|)=1}p)+1$.

\vspace{3mm}\noindent{\bf  The proof of Case 1:~~}%
Let $l=s+1$ and $r=(\prod_{i=1}^{i=s}|a_i|+k)$, where $k$ is the
least prime such that $(\prod_{i=1}^{i=s}|a_i|,k)=1$.

Set $e=\varphi (\prod_{p\leq
s^{l+1}\prod_{i=1}^{i=s}|a_i|+1,(p,\prod_{i=1}^{i=s}|a_i|)=1}p)+1$
and $c_i=r^{e^i}$ for $i=1,...,s+1$. By the known condition, Lemma 1
and Lemma 2, we know that there are positive integers $u, v$ such
that
$(\prod_{i=1}^{i=s}(a_i+b_ic_1u),\prod_{i=1}^{i=s}(a_i+b_ic_1v))=1$
and any prime divisor of
$\prod_{i=1}^{i=s}(a_i+b_ic_1u)\prod_{i=1}^{i=s}(a_i+b_ic_1v)$ is
greater than $a+b\times c_{s+1}$, where $a=\max\{|a_1|,...,|a_s|\},
b=\max\{b_1,...,b_s\}$. Let $P$ be the set of all prime divisors of
$\prod_{i=1}^{i=s}(a_i+b_ic_1u)$. And let $Q$ be the set of all
prime divisors of $\prod_{i=1}^{i=s}(a_i+b_ic_1v)$. We claim that
for every integer $m\geq
c=\prod_{i=1}^{i=s}(a_i+b_ic_1u)\prod_{i=1}^{i=s}(a_i+b_ic_1v)$,
there exists a positive integer $y$ such that
$f_1(y)>1,...,f_s(y)>1$ are all in $Z_m^*$. Certainly, in order to
prove that Case 1 holds, it is enough to consider the case $m=pqt$,
where $p\in P$, $q\in Q$ and $t\in N$. By the known condition and
Lemma 3, there must be a positive integer $x\leq p\times t$ such
that $(\prod_{j=1}^{j=s+1}\prod_{i=1}^{i=s}(a_i+c_j\times
b_ix),p\times t)=1$.  If $q|x$, then one can choose $y=c_1\times x$
and Case 1 holds. When $(q,x)=1$, let's consider the matrix
$$M=(m_{i,j})=\left(
\begin{matrix}
a_1+b_1c_1x, & \cdots,  & a_1+b_1c_{s+1}x\\
\cdots, & \cdots, & \cdots \\
a_s+b_sc_1x, & \cdots, & a_s+b_sc_{s+1}x
\end{matrix} \right).$$

Since $q>a+b\times c_{s+1}$, hence there is at most a number which
can be divided $q$ in each row of the matrix $M$. But there are
$s(s+1)$ elements in $M$. So, there must be some $j$ with $1\leq
j\leq s+1$ such that $(\prod_{i=1}^{i=s}(a_i+c_j\times b_ix), q)=1$.
This completes the proof of Case 1.

\vspace{3mm}\noindent{\bf  The proof of Case 2:~~}%
By the known conditions, it is easy to prove that there are positive
integers $u>1$ and $v>1$ such that for every $2\leq j\leq s$,
$a_j+b_ju>1$, and $( u\times\prod_{i=2}^{i=s}(a_i+b_iu),v\times
\prod_{i=2}^{i=s}(a_i+b_iv))=1$, and any prime divisor of $v\times
\prod_{i=2}^{i=s}(a_i+b_iv)$ is greater than $a+bsp\times
\prod_{i=2}^{i=s}|a_i|$, where $a=\max\{|a_1|,...,|a_s|\},
b=\max\{b_1,...,b_s\}$ and $p$ is the largest prime divisor of
$u\times \prod_{i=2}^{i=s}(a_i+b_iu)$. Let $P$ be the set of all
prime divisors of $u\times\prod_{i=2}^{i=s}(a_i+b_iu)$. And let $Q$
be the set of all prime divisors of $v\times
\prod_{i=2}^{i=s}(a_i+b_iv)$. We claim that for every integer $m\geq
c=u\times \prod_{i=2}^{i=s}(a_i+b_iu)\times
v\times\prod_{i=2}^{i=s}(a_i+b_iv)$, there exists a positive integer
$y$ such that $f_1(y)>1,...,f_s(y)>1$ are all in $Z_m^*$. Clearly,
in order to prove that it holds, it is enough to consider the case
$m=pqt$, where $p\in P$, $q\in Q$ and $t\in N$. By the known
conditions, there must be a positive integer $x\leq p\times t$ such
that $(x\times\prod_{i=2}^{i=s}(a_i+b_ix),p\times t)=1$.  Let
$z_k=kpt\prod_{i=2}^{i=s}|a_i|+x$ for $1\leq k\leq s+1$, and let's
consider the matrix
$$M=(m_{i,j})=\left(
\begin{matrix}
z_1, & \cdots,  & z_{s+1}\\
a_2+b_2z_1, & \cdots,  & a_2+b_2z_{s+1}\\
\cdots, & \cdots, & \cdots \\
a_s+b_sz_1, & \cdots, & a_s+b_sz_{s+1}
\end{matrix} \right).$$

Since $q>a+bsp\times \prod_{i=2}^{i=s}|a_i|$, hence there is at most
a number which can be divided $q$ in each row of the matrix $M$. But
there are $s(s+1)$ elements in $M$. So, there must be some $j$ with
$1\leq j\leq s+1$ such that $(z_j\times
\prod_{i=2}^{i=s}(a_i+z_j\times b_i), q)=1$. Let $y=z_j$ and this
shows that Case 2 holds and we completed the proof of equivalent
form of Dickson's conjecture.

\vspace{3mm}\noindent{\bf Remark 2:~~}%
In the aforementioned proof on equivalent form of Dickson's
conjecture, we only consider the existence of constant $c$ and give
a rough estimation. So, our proof is not good and there must be
simpler proofs. The problem of the low bound of $c$ is very
interesting. In [15], we have discussed this problem. Unfortunately,
this problem is difficult and we only obtain several results in some
special cases. For example, let $L_{f_1,...,f_s}(c)$ be the least
positive integer in the equivalent form of Dickson's conjecture. We
have $L_{f_1,f_2}(c)=7$ when $f_1=x$ and $f_2=x+2$.
$L_{f_1,f_2}(c)=16$ when $f_1=x$ and $f_2=2x+1$. For the details,
see [16-17]. By using the methods in this paper, one can give
simpler proofs and get slightly stronger results as follows:

\vspace{3mm}\noindent{\bf Corollary 1:~~}%
Let $f_1=x$ and $f_2=x+2$. For every positive integer $m>10$, there
is a positive integer $x\equiv 5(\mod 6)$ such that $x\in Z_m^*$ and
$x+2\in Z_m^*$.

\vspace{3mm}\noindent{\bf The proof of Corollary 1:~~}%
When $20>m>10$, it is clear. Let $m\geq 20$. If $(m,35)=1$, then we
choose $x=5$ and if  $(m,323)=1$, then we choose $x=17$. So, it is
enough to consider the case $m=pqt$, where $p$ is a prime in
$\{5,7\}$ and $q$ is a prime in $\{17,19\}$, and $t$ is a positive
integer. Note that $6pt-1$,$6pt+1$,$12pt-1$ and $12pt+1$ are
pairwise relatively prime. Therefore, either both of $6pt-1$ and
$6pt+1$ are co-prime to $m=pqt$, or both of $12pt-1$ and $12pt+1$
are co-prime to $m=pqt$. Note that
$1<6pt-1<6pt+1<12pt-1<12pt+1<17pt\leq pqt$. Thus Corollary 1 holds.

\vspace{3mm}\noindent{\bf Corollary 2:~~}%
Let $f_1=x$ and $f_2=2x+1$. For every positive integer $m>45$, there
is a positive integer $x\equiv 5(\mod 6)$ such that $x\in Z_m^*$ and
$2x+1\in Z_m^*$.

\vspace{3mm}\noindent{\bf The proof of Corollary 2:~~}%
When $168>m>45$, one can  prove directly. Let $m\geq 168$. It is
enough to consider the case $m=pqt$, where $p$ is a prime in
$\{5,11\}$ and $q$ is a prime in $\{83,167\}$, and $t$ is a positive
integer. Note that $6pt-1$,$12pt-1$,$30pt-1$ and $60pt-1$ are
pairwise relatively prime. Therefore, either both of $6pt-1$ and
$12pt-1$ are co-prime to $m=pqt$, or both of $30pt-1$ and $60pt-1$
are co-prime to $m=pqt$. Noth that
$1<6pt-1<12pt-1<30pt-1<60pt-1<83pt\leq pqt$. Thus Corollary 2 holds.

\vspace{3mm}\noindent{\bf Remark 3:~~}%
For any given positive integer $k$, let $f_1=x$ and $f_2=x+2k$, we
refined the result in [16] without using the method in this paper.
But we did not obtain the precise value of $L_{f_1,f_2}(c)$. We will
continue to consider this problem.

\vspace{3mm}\noindent{\bf Remark 4:~~}%
Dickson's conjecture implies that for any $n>4$,
$f_1(x)=1+2x$,$f_2(x)=1+2^2x$,$f_3(x)=1+2^{2^2}x$,
$f_4(x)=1+2^{2^3}x$,$f_5(x)=1+2^{2^4}x$,...,$f_{n+1}(x)=1+2^{2^n}x$
can represent infinitely many prime for $x$. Naturally, we want to
know the least $x$ such that $f_1(x),...,f_{n+1}(x)$ represent
simultaneously primes. This leads to another interesting problem---
estimating the bound of the least $x$.

\vspace{3mm}\noindent{\bf Remark 5:~~}%
Dickson's conjecture implies the first Hardy-Littlewood  conjecture
[9] which states that for $s\geq 1$, $f_i(x)=1+2b_ix$ for
$i=1,...,s$. If there does not exist any integer $n>1$ dividing all
the products $k\prod_{i=1}^{i=s}f_i(k)$, for every integer $k$, then
the number of primes $p\leq x$ such that $p+2b_1,...,p+2b_s$  are
all prime is about $\int_2^x\frac{dt}{(\log t)^{s+1}}$. The second
Hardy-Littlewood conjecture states that $\pi (x+y)\leq \pi (x)+\pi
(y)$ for all $x,y\geq 2$, where $\pi (x)$ is the prime counting
function. It's interesting that in 1974, Richards [18] proved that
the first and second conjectures are incompatible with each other.
Generally speaking, there is  an increased tendency to disprove the
second Hardy-Littlewood conjecture. In 1962, Paul T. Bateman and
Roger A. Horn gave a heuristic asymptotic formula concerning the
distribution of prime numbers, which implies the first
Hardy-Littlewood conjecture. What is more interesting is of that
Friedlander John and Granville Andrew [20-22] showed that
Bateman-Horn's asymptotic formula does not always hold. In 2006, by
considering the case of system of non-constant affine-linear forms,
Green and Tao [13] generalized the first Hardy-Littlewood
conjecture. Based on their work, the author got a heuristic
conclusion that Dickson's conjecture must hold and it can be
generalized to the multivariable case. Next section, we will try to
give the precise Dickson's conjecture in the multivariable case.

\section{The generalization of Dickson's conjecture}
Let $s,k\geq 1$ be positive integers and let
$f_1(x_1,...,x_k),...,f_s(x_1,...,x_k)$ be number-theoretic
functions from $N^k$ to $Z$. To begin with, we explain the meaning
that $f_1(x_1,...,x_k),...,f_s(x_1,...,x_k)$  represent
simultaneously primes for infinitely many integral points
$(x_1,...,x_k)$. We must point out that we do not consider the
trivial cases such as $f_1(x)=\left\{
\begin{array}{c}
2,(x,2)=1 \\
1+7x,2|x\\
\end{array}
\right. $ and $f_2(x)=\left\{
\begin{array}{c}
6x+1,(x,2)=1 \\
5,2|x\\
\end{array}
\right. $, and so on. Clearly, $f_1(x), f_2(x)$ represent
simultaneously primes for infinitely many $x$. But, this is not the
case considered. Considering that
$f_1(x_1,...,x_k),...,f_s(x_1,...,x_k)$ represent simultaneously
primes for infinitely many integral points $(x_1,...,x_k)$, we
firstly hope that for any $1\leq i\leq s$, $f_i(x_1,...,x_k)$ itself
can represent primes for infinitely many integral points
$(x_1,...,x_k)$. So, we do not consider the case that for some
$1\leq i\leq s$, $f_i(x_1,...,x_k)$ is a constant. Moreover, we
require that there must be an infinite sequence of integral points
$(x_{11},...,x_{k1})$, ..., $(x_{1i},...,x_{ki})$, ... such that for
any $i\neq j$, $f_1(x_{1i},...,x_{ki})\neq
f_1(x_{1j},...,x_{kj}),...,f_s(x_{1i},...,x_{ki})\neq
f_s(x_{1j},...,x_{kj})$ hold simultaneously.

\vspace{3mm}Therefore, the meaning that $f_1(x_1,...,x_k)$, ...,
$f_s(x_1,...,x_k)$  represent simultaneously primes for infinitely
many integral points $(x_1,...,x_k)$ is of that
$f_1(x_1,...,x_k),...,f_s(x_1,...,x_k)$ represent distinct primes
for infinitely many integral points $(x_1,...,x_k)$ respectively,
moreover, there is an infinite sequence of integral points
$(x_{11},...,x_{k1})$, ..., $(x_{1i},...,x_{ki})$, ... such that for
any positive integer $r$, $f_1(x_{1r},...,x_{kr})$,...,
$f_s(x_{1r},...,x_{kr})$ represent simultaneously primes, and for
any $i\neq j$, $f_1(x_{1i},...,x_{ki})\neq f_1(x_{1j},...,x_{kj})$,
..., $f_s(x_{1i},...,x_{ki})\neq f_s(x_{1j},...,x_{kj})$ hold
simultaneously.

\vspace{3mm} Thus, if $f_1(x_1,...,x_k)$, ..., $f_s(x_1,...,x_k)$
represent simultaneously primes for infinitely many integral points
$(x_1,...,x_k)$, then there is an infinite sequence of integral
points $(x_{11},...,x_{k1})$, ..., $(x_{1i},...,x_{ki})$, ... such
that for any positive integer $r$, $f_1(x_{1r},...,x_{kr})$,...,
$f_s(x_{1r},...,x_{kr})$ represent simultaneously primes, and for
any positive integer $c$, there is a positive integer $l$ such that
for every integer $m\geq l$, we have $f_1(x_{1m},...,x_{km})\geq
c$,..., $f_s(x_{1m},...,x_{km})\geq c$ and
$f_1(x_{1m},...,x_{km})$,..., $f_s(x_{1m},...,x_{km})$ represent
simultaneously primes.

\vspace{3mm} So,  we get a natural necessary condition that
$f_1(x_1,...,x_k)$, ..., $f_s(x_1,...,x_k)$ represent simultaneously
primes for infinitely many integral points $(x_1,...,x_k)$ as
follows: there exists an infinite sequence of integral points
$(x_{11},...,x_{k1})$, ..., $(x_{1i},...,x_{ki})$, ... such that
$\prod_{j=1}^{j=s}f_j(x_{11},...,x_{k1})$ , ... ,
 $\prod_{j=1}^{j=s}f_j(x_{11},...,x_{ki})$, ... are pairwise
relatively prime and $f_j(x_{11},...,x_{ki})>1$ for each $i$ and
$j$. In [15], we refined this necessary condition and obtained the
following conjectures 2 and 3:

\vspace{3mm}\noindent{\bf Conjecture 2:~~}%
Let $f_1(x_1,...,x_k),...,f_s(x_1,...,x_k)$  be $s$ multivariable
number-theoretic functions. If $f_1(x_1,...,x_k)$, ...,
$f_s(x_1,...,x_k)$ represent simultaneously primes for infinitely
many integral points $(x_1,...,x_k)$, then there is always a
constant $c$  such that for every positive integer $m>c$, there
exists an integral point $(y_1,...,y_k)$ such that
$f_1(y_1,...,y_k)>1,...,f_s(y_1,...,y_k)>1$ are all in $Z_m^* $.

\vspace{3mm} In Section 1 of this paper, we have proved Conjecture 2
holds when $ k=1$ and $f_1(x_1),...,f_s(x_1)$ are linear
polynomials.

\vspace{3mm}\noindent{\bf  Conjecture 3:~~}%
Let $f_1(x_1,...,x_k),...,f_s(x_1,...,x_k)$ be multivariable
polynomials with  integral coefficients. If there is a positive
integer $c$ such that for every positive integer $m\geq c$, there
exists an integral point $(y_1,...,y_k)$ such that
$f_1(y_1,...,y_k)>1,...,f_s(y_1,...,y_k)>1$ are all in $Z_m^* $, and
there exists an integral point $(z_1,...,z_k)$ such that
$f_1(z_1,...,z_k)\geq c,...,f_s(z_1,...,z_k)\geq c$ are all primes,
then $f_1(x_1,...,x_k)$, ..., $f_s(x_1,...,x_k)$  represent
simultaneously primes for infinitely many integral points
$(x_1,...,x_k)$.

\vspace{3mm} As a special case of Conjecture 3, in 1997, Fouvry,
Etienne and Iwaniec, Henryk [28] proved that
$f_1(x_1,x_2)=x_1^2+x_2^2$ and $f_2(x_1,x_2)=1\times x_1+0\times
x_2$ represent simultaneously primes for infinitely many integral
points $(x_1,x_2)$.

\vspace{3mm} Now, let $f_1(x_1,...,x_k),...,f_s(x_1,...,x_k)$ be
multivariable polynomials of degree 1 with  integral coefficients.
If  there is a positive integer $c$ such that for every positive
integer $m\geq c$, there exists an integral point $(y_1,...,y_k)$
such that $f_1(y_1,...,y_k)>1,...,f_s(y_1,...,y_k)>1$ are all in
$Z_m^* $, then for any $1\leq i\leq s$, the greatest common divisor
of coefficients of $f_i(x_1,...,x_k)$ is 1. So  $f_i(x_1,...,x_k)$
represents  primes for infinitely many integral points
$(x_1,...,x_k)$. Of course, there exists an integral point
$(z_1,...,z_k)$ such that $f_i(z_1,...,z_k)\geq c$ is prime. We
conjecture that there exists an integral point $(w_1,...,w_k)$ such
that $f_1(w_1,...,w_k)\geq c,...,f_s(w_1,...,w_k)\geq c$ are all
primes. By Conjecture 3, we get a generalization of Dickson's
conjecture as follows:

\vspace{3mm}\noindent {\bf The generalization of Dickson's conjecture:~~}%
Let $s, k \in N$ and let $f_1(x_1,...,x_k),...,f_s(x_1,...,x_k)$  be
multivariable polynomials of degree 1 with  integral coefficients,
if there is a positive integer $c$ such that for every positive
integer $m\geq c$, there exists an integral point $(y_1,...,y_k)$
such that $f_1(y_1,...,y_k)>1,...,f_s(y_1,...,y_k)>1$ are all in
$Z_m^*$, then $f_1, ... , f_s$ represent simultaneously primes for
infinitely many integral points $(x_1,...,x_k)$.

\vspace{3mm} Let's consider several simple cases of linear system
$L$ which satisfies the condition of the generalization of Dickson's
conjecture, where
$$L=\left\{
\begin{array}{c}
f_1(x_1,...,x_k)=a_{11}x_1+...+a_{1k}x_k +b_1\\
...........................................................\\
f_s(x_1,...,x_k)=a_{s1}x_1+...+a_{sk}x_k +b_s\\
\end{array}
\right. .$$

\vspace{3mm} For example: let $$L_1=\left\{
\begin{array}{c}
f_1(x_1,x_2)=x_1+2x_2 \\
f_2(x_1,x_2)=2x_1+x_2 \\
\end{array}
\right. .$$

By Corollary 1, one can prove that there is a positive integer $c=6$
such that for every positive integer $m> c$, there exists an
integral point $(y_1,y_2)$ such that $f_1(y_1,y_2)>1,f_2(y_1,y_2)>1$
are all in $Z_m^*$. By the generalization of Dickson's conjecture,
$f_1(x_1,x_2)=x_1+2x_2, f_2(x_1,x_2)=2x_1+x_2$ should represent
simultaneously primes for infinitely many integral points
$(x_1,x_2)$. Clearly, this is a necessary condition that there are
infinitely many twin primes.

As another example, let's consider $$L_2=\left\{
\begin{array}{c}
f_1(x_1,x_2)=x_1+x_2+1 \\
f_2(x_1,x_2)=3x_1+x_2+3 \\
\end{array}
\right. .$$

By Corollary 2, one can prove that there is a positive integer $c=7$
such that for every positive integer $m> c$, there exists an
integral point $(y_1,y_2)$ such that $f_1(y_1,y_2)>1,f_2(y_1,y_2)>1$
are all in $Z_m^*$. By the generalization of Dickson's conjecture,
$f_1(x_1,x_2)=x_1+x_2+1, f_2(x_1,x_2)=3x_1+x_2+3$ should represent
simultaneously primes for infinitely many integral points
$(x_1,x_2)$. Clearly, this is a necessary condition that there are
infinitely many safe primes or Sophie-Germain primes.

\vspace{3mm}\noindent{\bf Remark 6:~~}%
We are very interesting in the following linear system:
$$L_3=\left\{
\begin{array}{c}
f_1(x_1,...,x_k)=a_{11}x_1+...+a_{1k}x_k\\
.....................................................\\
f_k(x_1,...,x_k)=a_{k1}x_1+...+a_{kk}x_k\\
\end{array}
\right. .$$

We write $L_3=AX^T$, where $X=(x_1,...,x_k)^T,A=\left(
\begin{matrix}
a_{11}, & \cdots,  & a_{1k}\\
\cdots, & \cdots, & \cdots \\
a_{k1}, & \cdots, & a_{kk}
\end{matrix} \right)$ with $a_{ij}\in Z$.

Clearly, the matrix $A$ decides whether $f_1(X), ... , f_k(X)$ in
the system $L_3$ represent simultaneously primes. When
$f_1(X)=p_1,...,f_k(X)=p_k$ in the system $L_3$ represent
simultaneously primes, we also say that $A$ represents a prime
vector $(p_1,..., p_k)$. If there is a positive integer $c$ such
that for every positive integer $m> c$, there exists an integral
point $Y$ such that $f_1(Y)>1,...,f_s(Y)>1$ are all in $Z_m^*$, we
say that $A$ has the good property. The generalization of Dickson's
conjecture states that a sufficient and necessary condition that $A$
represents infinitely many prime vectors is of that $A$ has the good
property. Surely, a unit matrix has the good property, and all
matrixes which have the good property do not always form a group.
However, there must be many interesting properties on $A$. These
problems are left as interesting exercises.

\vspace{3mm}\noindent{\bf Remark 7:~~}%
By the idea in Section 9 of [15], one can get another sufficient and
necessary condition that $A$ represents infinitely many prime
vectors is of that there is a positive integer $c$ such that for
every positive integer $m\geq c$, there exists an integral point
$(y_1,...,y_k)$ such that
$f_1(y_1,...,y_k)>1,...,f_s(y_1,...,y_k)>1$ are all in $Z_{m!}^* $.
This sufficient and necessary condition implies that there exists an
integral point $(x_1,...,x_k)$ such that all of
$f_1(x_1,...,x_k),...,f_s(x_1,...,x_k)$  in $Z_{r!}^* $ are prime,
where $r$ satisfies $(r-1)!<L_{f_1,...,f_s}(c)\leq r!$, and
$L_{f_1,...,f_s}(c)$ is the least value of $c$. It shows that the
problem of determining the low bound of constant $c$ is of critical
importance, and it perhaps leads to a new way for proving Dickson's
conjecture and its generalization.

\vspace{3mm}\noindent{\bf Remark 8:~~}%
Finally,  from a computational point of view, the author asks two
questions to close this section:

Q 1: Let $s, k \in N$ and let
$f_1(x_1,...,x_k),...,f_s(x_1,...,x_k)$ be multivariable polynomials
of degree 1 with  integral coefficients. Is there an efficient
algorithm for determining whether there is a positive integer $c$
such that for every positive integer $m\geq c$, there exists an
integral point $(y_1,...,y_k)$ such that
$f_1(y_1,...,y_k)>1,...,f_s(y_1,...,y_k)>1$ are all in $Z_m^* $?
Especially, for any given the matrix $A$, is there an efficient
algorithm for determining whether $A$ has the good property?

Q 2: For any given the matrix $A$ which has the good property, is
there an efficient algorithm for finding $L_{f_1,...,f_s}(c)$?

\section{Conclusion}
This paper should be a part of the paper \emph{On the infinitude of
some special kinds of primes}. But personally I'd prefer to let it
become an independent paper. From a computational point of view,
what mankind can solve ideally is only a linear system for thousands
of years. In the non-linear case, there is not a generic method.
Even solving a system of multivariate equations of order 2 over the
finite field $GF(2)$ is NP-hard [23]. This maybe is God's punishment
for mankind. In the author's eyes, the problems in nonlinear systems
are rather unattackable. Therefore, the author prefers considering
some question in the linear systems to studying those in the
nonlinear systems although people have made a lot of great progress.
At least, the author feels that to solve twin prime conjecture
perhaps is easier than to solve Landau's first conjecture. Although
it has been proved that $x^2+y^4$ and $x^3+2y^3$ can represent
infinitely many primes respectively [24-25], it perhaps is
"unattackable at the present state of science" to prove that
$x^2+y^4$ and $x^3+2y^3$ can simultaneously represent infinitely
many primes. And it is only the personal viewpoint. Number theorists
will tell us the correct answer.

\vspace{3mm} Of course, in the linear case, there are some very hard
problems, too. For instance, by Euclid's algorithm, one can solve
the following problem: for any given integers $a_1,...a_n$, find
integers $x_1,...x_n$ with polynomial time such that
$a_1x_1+...+a_nx_n=(a_1,...a_n)$. But, when $n$ is large, there is
not an efficient algorithm for finding $x_1,...x_n$ such that
$a_1x_1+...+a_nx_n=(a_1,...a_n)$ and
$((x_1)^2+...+(x_n)^2)^{\frac{1}{2}}$ is the least. This problem
also is NP-hard. Moreover, Majewski and Havas [26] proved that
finding the minimum $\max_{1\leq i\leq n}|x_i|$ is NP-complete. As
another example, Hilbert [27] asked whether the diophantine equation
$ax+by+c=0$ is always solvable in primes $x,y$ if $a,b,c$ are given
pairwise relatively prime integers. Hilbert's problem is still open
in spite of many excellent mathematicians have been made unremitting
endeavor. We view these puzzles as a special kind of linear
problems. Although the author does not know whether Dickson's
conjecture and its generalization belongs this special kind
problems, he believes that they must be solved perfectly in the near
future.

\section{Acknowledgements}
 Erd\"{o}s said "the meaning of life is to prove and conjecture". Before being a father,
 I viewed his words as Holy Bible. At my QQ, on my wall, in my mind,
 I always wrote them in order to urge myself. However, after being a father, I did not think so.  For me,
Erd\"{o}s was only half right. In my eyes, the other half meaning of
life is to honor father and educate son. Confucius said "the father
who does not teach his son his duties is equally guilty with the son
who neglects them." Therefore, we must educate our sons. As a son,
we should honor our parents. "For thousands of years, father and son
have stretched wistful hands across the canyon of time." Alan
Valentine said.

\vspace{3mm} George Herbert said "one father is more than a hundred
schoolmasters." Thank my father Hanchao Zhang for his deep love.
Thank him  for teaching me how to live and for encouraging me to
study, especially, for teaching me how to persistently think  in my
research area. Today is Father's Day. I am very sorry to say that I
haven't seen him for a long time. Now my wish is to go home as
quickly as possible and look my father up. Hope that he has a good
health!

\section{Appendix}
Last Sunday was my father's birthday. I am very sorry to say that I
didn't go home to see him again. I finished writing the draft
hurriedly to celebrate his birthday. Below is that self-contained
paper which has no theorems but dreams.

\vspace{3mm}\noindent{\textbf{ A brief introduction to the theory
that several multivariable integral polynomials simultaneously
represent infinitely many primes\\-------------------{\small\it
Dedicated to my father on the occasion of his 66th birthday}}

\vspace{3mm}\noindent{\textbf{Abstract}}

\vspace{3mm}In 2006,  by generalizing Hardy-Littlewood's formula,
Green and Tao were the first to consider the question that a system
of non-constant affine-linear forms from $Z^n$ to $Z^m$ represents
infinitely many primes. In this appendix, we further generalize
their conjecture and introduce briefly the basic theory that several
multivariable integral polynomials represent simultaneously prime
numbers for infinitely many integral points.

\vspace{3mm}\noindent{\textbf{Some basic notations and definitions}}

\vspace{3mm}Let $Z$ be the set of integers. Denote the set of all
prime numbers by $P$ and denote the set of all natural numbers or
positive integers by $N$. We define affine $n$-space over $Z$,
denoted $Z^n$, to be the set of all $n$-tuples of elements of $Z$.
An element $x=(x_1,...,x_n)$ in $Z^n$ is called an integral point
and $x_i$ is called the coordinates of $x$. Let $Z[x_1,...,x_n]$ be
the polynomial ring in $n$ variables over $Z$. Let $Z_n^*=\{x\in
N|1\leq x\leq n, \gcd(x,n)=1\}$ be the set of positive integers less
than or equal to $n$ that are coprime to $n$. Let $Z_n=\{x\in
Z|0\leq x\leq n-1\}$.

\vspace{3mm} Let's consider the map $F: Z^n\rightarrow Z^m$ for all
integral points $x=(x_1,...,x_n)\in Z^n$, $F(x)=(f_1(x),...,f_m(x))$
for distinct polynomials $f_1,...,f_m\in Z[x_1,...,x_n]$.  In this
case, we call $F$ a polynomial map on $Z^n$. Call $F$ a polynomial
map on $N^n$ if $F: N^n\rightarrow Z^m$ for all integral points
$x\in N^n$, $F(x)=(f_1(x),...,f_m(x))$.  We call $F$ (on $Z^n$ or
$N^n$ ) linear or a system of non-constant affine-linear forms if
for any $1\leq i\leq m$, $f_i(x_1,...,x_n)$ has degree 1. We call
the polynomial map $F$ on $Z^n$ ( or $N^n$ ) admissible if for every
positive integer $r$ there exists an integral point
$x=(x_1,...,x_n)\in Z^n$ ( or $N^n$ ) such that $r$ is coprime to
$f_1(x)\times...\times f_m(x)=\prod_{i=1}^{i=m}f_i(x)$, moreover,
$f_i(x)>1$ for $1\leq i\leq m$. We call the polynomial map $F$ on
$Z^n$ ( or $N^n$ ) strongly admissible if there is a positive
integer $C$ such that for every positive integer $k\geq C$, there
exists an integral point $x=(x_1,...,x_n)$ such that
$f_1(x)>1,...,f_m(x)>1$ are all in $Z_k^*$. We call the least
positive integer $C$ such that $F$ is strongly admissible a strongly
admissible constant.

\vspace{3mm} Let $S$ be set of all solutions of the simultaneous
equations: $f_1(x_1,...,x_n)\in P,...,f_m(x_1,...,x_n)\in P$. Let
$H=F(S)\in P^m$ be the image of $S$ under $F$. An element of $H$ is
called a prime point. Let
$\Omega_F(\alpha)=\#\{(f_1(x),...,f_m(x))\in P^m:f_i(x)\leq \alpha
(1\leq i\leq m)\}$ be the number of prime points
$(f_1(x),...,f_m(x))$ whose any coordinate $f_i(x)$ is less than or
equal to a positive real number $\alpha$. Let
$\Psi_F(\beta)=\#\{x=(x_1,...,x_n)\in Z^n : (f_1(x),...,f_m(x))\in
P^m, |x_i|\leq \beta (1\leq i\leq n)\}$ be the number of integral
points $x=(x_1,...,x_n)$ such that the  absolute value of any
coordinate $|x_i|$ is less than or equal to a non-negative real
number $\beta$ such that $(f_1(x),...,f_m(x))\in P^m$. More
generally,  one would like to consider
$\Omega_F(\alpha_1,...,\alpha_m)=\#\{(f_1(x),...,f_m(x))\in
P^m:f_i(x)\leq \alpha_i (1\leq i\leq m)\}$ for $m$ positive real
numbers $\alpha_1,...,\alpha_m$ and
$\Psi_F(\beta_1,...,\beta_n)=\#\{x=(x_1,...,x_n)\in Z^n :
(f_1(x),...,f_m(x))\in P^m, |x_i|\leq \beta_i (1\leq i\leq n)\}$ for
$n$ non-negative real numbers $\beta_1,...,\beta_n$.

\vspace{3mm}\noindent{\textbf{Two basic questions}}

\vspace{3mm}Let's begin with two well-known results in Analytic
Number Theory. In 1988,  Friedlander and Iwaniec proved that the
polynomial $f_1(x,y)=x^2+y^4$ represents infinitely many primes.
Three years later Heath-Brown showed that $f_2(x,y)=x^3+2y^3$ also
represents infinitely many primes. Naturally, one might ask: do
$f_1(x,y)=x^2+y^4$ and $f_2(x,y)=x^3+2y^3$ represent simultaneously
primes for infinitely many integral points $(x,y)\in N^2$? More
generally, for a given polynomial map $F$ on $Z^n$, how to determine
whether $H$ is an infinite set or not? Concerning the theory that
several multivariable integral polynomials represent simultaneously
prime numbers for infinitely many integral points, generally
speaking, the main goal is to consider the following two questions:

\vspace{3mm}{\bf Question 1:~~}%
Let $F$ be a given polynomial map on $Z^n$, how to determine whether
$H$ is an infinite set or not? Especially, how to determine whether
$H$ is empty?

\vspace{3mm}{\bf Question 2:~~}%
Let $F$ be a given polynomial map on $Z^n$. We assume that $H$ is
infinite. How to estimate $\Omega_F(\alpha)$,
$\Omega_F(\alpha_1,...,\alpha_m)$ (or $\Psi_F(\beta)$,
$\Psi_F(\beta_1,...,\beta_n)$) and give a good approximation to
$\Omega_F(\alpha)$, $\Omega_F(\alpha_1,...,\alpha_m)$ (or
$\Psi_F(\beta)$, $\Psi_F(\beta_1,...,\beta_n)$)?

\vspace{3mm} In this appendix, we consider mainly Question 1.
Incidentally we mention some  heuristic formulae on Question 2.

\vspace{3mm}\noindent{\textbf{Some heuristic conjectures on Question
1}}

\vspace{3mm} On Question 1, Green and Tao is the first to consider
the case that $F$ is a system of non-constant affine-linear forms on
$Z^n$ and generalize Hardy-Littlewood Conjecture. They believe that
the following conjecture 1 holds.

\vspace{3mm} {\bf Conjecture 1 (Green-Tao):~~}%
If a system of affine-linear forms $F$ on $Z^n$ is admissible, then
$H$ is infinite.

\vspace{3mm}In the case that $F$ is a system of non-constant
affine-linear forms on $N^n$, in 2009, we conjecture that if $F$ is
strongly admissible, then $H$ is infinite. This can be further
generalized as follows:

\vspace{3mm}{\bf Conjecture 2:~~}%
If a system of affine-linear forms $F$ on $Z^n$ is strongly
admissible, then $H$ is infinite.

\vspace{3mm}In our subsequent paper \emph{An equivalent form of
Green-Tao Conjecture}, we will prove that Conjecture 1 and
Conjecture 2 are equivalent using Generalized Chinese Remainder
Theorem which states that if $\gcd (a,b)=1$, the Cartesian product
$Z_a^*\times Z_b^*$ ($Z_a\times Z_b$) is isomorphic to $Z_{ab}^*$
($Z_{ab}$). More generally, if $\gcd (a_i,a_j)=1$ for $1\leq i\neq
j\leq n$, then $(Z_{a_1}^*)^m\times ...\times (Z_{a_n}^*)^m\simeq
(Z_{a_1...a_n}^*)^m$  or $(Z_{a_1})^m\times ...\times
(Z_{a_n})^m\simeq (Z_{a_1...a_n})^m$ for any positive integer $m$.

\vspace{3mm} We believe that the  expression of Conjecture 2 is
better than the expression of Conjecture 1 because it is more
convenient to generalize by using the notation "strongly
admissible". In fact, one could generalize Conjecture 2 to the case
that $F$ is a polynomial map on $Z^n$ as follows.

\vspace{3mm}{\bf Conjecture 3 (The generalization of H Hypothesis):~~}%
Let $F$ be a polynomial map on $Z^n$. If $F$ is strongly admissible,
and there exists an integral point $y=(y_1,...,y_n)$ such that
$f_1(y)\geq C,..., f_m(y)\geq C$ are all primes, then $H$ is
infinite, where $C$ is the strongly admissible constant.

\vspace{3mm} Conjecture 3 actually gives the sufficient and
necessary condition that $H$ is infinite for the given polynomial
map $F$ on $Z^n$.

\vspace{3mm} Historically, on Question 1, it goes back to
Bouniakowsky even Dirichlet. In 1857, Bouniakowsky conjectured the
following:

\vspace{3mm}{\bf Conjecture 4 (Bouniakowsky):~~}%
If $f(x)$ is an irreducible polynomial with integral coefficients
and positive leading term and $f(x)$ is admissible, then $f(x)$ is
prime for an infinite number of integers $x$.

\vspace{3mm} Bouniakowsky perhaps is the first person who finds that
the condition "admissible" is necessary.

\vspace{3mm}Concerning the simultaneous values of several linear
polynomials, Dickson stated a very important  conjecture in 1904. In
1958,  Schinzel and Sierpinski further generalized Bouniakowsky's
conjecture and Dickson's conjecture as follows:

\vspace{3mm}{\bf Conjecture 5 (H Hypothesis):~~}%
Let $f_1(x),...,f_m(x)$ be irreducible polynomials with integral
coefficients and positive leading coefficient. If the polynomial map
$F:N\rightarrow Z^m$ for all $x\in N$, $F(x)=(f_1(x),...,f_m(x))$ is
admissible, then there exist infinitely many natural numbers $r$
such that all numbers $f_1(r),...,f_m(r)$ are primes.

\vspace{3mm} Conjecture 5 is the special case of Conjecture 3. Thus,
Conjecture 5 can be re-stated: Let $F:N\rightarrow Z^m$ be a
polynomial map. If $F$ is strongly admissible, and there exists an
integral point $y=(y_1,...,y_n)$ such that $f_1(y)\geq C,...,
f_m(y)\geq C$ are all primes, then $H$ is infinite, where $C$ is the
strongly admissible constant.

\vspace{3mm} Conjecture 5 implies many brilliant results. For
example, the primes contain arbitrarily long arithmetic progressions
(or polynomial progressions) because the polynomial map
$F:N\rightarrow Z^m$ for all $x\in N$, $F(x)=(x+m!,...,x+m\times
m!)$ is admissible (or $F:N\rightarrow Z^m$ for all $x\in N$,
$F(x)=(xg_1(x)+(m!)^k+1,...,xg_m(x)+(m!)^k+1)$ is admissible.).

\vspace{3mm}\noindent{\textbf{Some heuristic formulae on Question
2}}

\vspace{3mm} On  Question 2, it goes back to Hardy and Littlewood
even Legendre and Gauss. A variant of Hardy-Littlewood's conjecture
is the following:

\vspace{3mm}{\bf Conjecture 6 (Hardy-Littlewood):~~}%
Let $f_1(x)=a_1x+b_1,...,f_m(x)=a_mx+b_m$ be linear polynomials with
integral coefficients and positive leading coefficient. If the
polynomial map $F:N\rightarrow Z^m$ for all $x\in N$,
$F(x)=(f_1(x),...,f_m(x))$ is admissible, then $\Psi_F(\beta)\sim
\{\prod_p
(1-\frac{v(p)}{p})(1-\frac{1}{p})^{-m}\}\frac{\beta}{\log^m \beta}$,
where for each prime $p$, $v(p)$ is the number of distinct $k (\mod
p)$ for which $p$ does divide $\prod_{i=1}^{i=m} (a_ik+b_i)$.

\vspace{3mm}In 1962, Bateman and Horn generalized Conjecture 6 and
obtained the following conjecture 7.

\vspace{3mm}{\bf Conjecture 7 (Bateman-Horn):~~}%
Let $f_1(x),...,f_m(x)$ be irreducible polynomials with integral
coefficients and positive leading coefficients. If the polynomial
map $F:N\rightarrow Z^m$ for all $x\in N$,
$F(x)=(f_1(x),...,f_m(x))$ is admissible, then $\Psi_F(\beta)\sim
\frac{1}{\prod_{i=1}^{i=m}d_i}\{\prod_p
(1-\frac{v(p)}{p})(1-\frac{1}{p})^{-m}\}\frac{\beta}{\log^m \beta}$,
where $d_i =\deg f_i(x)$ and for each prime $p$, $v(p)$ is the
number of distinct $k (\mod p)$ for which $p$ does divide
$\prod_{i=1}^{i=m} f_i(k)$.

\vspace{3mm} However, Friedlander John and Granville Andrew showed
that Conjecture 7 does not always hold.

\vspace{3mm}In 2006, Green and Tao further generalized Conjecture 6
to the case of affine $n$-space over $Z$.  Thus, one will expect
that if Conjecture 3 holds, then there will be an approximation to
$\Omega_F(\alpha)$ (or $\Psi_F(\beta)$) which will further
generalize the work of Green and Tao. What does this asymptotic
formula look like? We left this question to the readers who are
interested in this area. However, it miles to go. Like the early
stage of some theories (such as Algebraic Geometry), at present it
is full of many heuristic conjectures in this area. Of course, we
believe that correct conjectures will come and the serious theory
will be establish systematically and scientifically in the future.
Let's try the item to wait.

\clearpage
\end{document}